\documentclass[12pt]{article}

\usepackage{graphicx}
\usepackage{amsmath}    
\usepackage{verbatim}   
\usepackage{amsthm}
\usepackage{amsfonts}
\usepackage{amssymb}
\usepackage{hyperref}   
\usepackage{listings}
\usepackage[linesnumbered,ruled,vlined,commentsnumbered]{algorithm2e}
\usepackage[export]{adjustbox}

\begin{document}
   \begin{center}
 
{\bf\large On Genus $g$ Orientable Crossing Numbers of Small Complete Graphs} \\
 
       \vspace{0.5cm}

       \textbf{Yoonah Lee}
              \vspace{0.5cm}

\end{center}

\begin{abstract}
The current state of knowledge of genus $g$ orientable crossing numbers of complete graphs through $K_{11}$ is reviewed. It is shown that $cr_3(K_{10})=3,$ $cr_3(K_{11})\leq14,$ and $cr_4(K_{11})=4.$ It is established with the aid of an algorithm that there are precisely two non-isomorphic embeddings of $K_9$ with a hexagon with all its vertices distinct on a surface of genus 3.
\end{abstract}

\section*{Introduction}
Given a surface $S$ and a graph $G$, a \emph{drawing} of $G$ has edges with two distinct vertices as endpoints and no vertices that share three edges \cite{szekely}. The \emph{crossing number of $G$ on $S$}, denoted $cr_S(G)$, is the minimum number of crossings in any drawing of $G$ on $S$. \emph{Good drawings} are considered to have three conditions: no arc connecting two nodes intersects itself, two arcs have at most one point of intersection, and two arcs with the same node have no other intersecting points. \emph{Optimal drawings} have the minimum number of crossings of $G$ on $S$ \cite{guy1971}. It can be shown that optimal drawings are good drawings.

The crossing numbers of complete graphs $K_n$ on orientable surfaces $S_g$ of positive genus has been an object of interest since Heawood's map-coloring conjecture in 1890 \cite{heawood}. Despite significant progress, the values of the crossing number for even rather small complete graphs on these surfaces has been unclear. But when $K_{n-1}$ can be embedded in $S_g,$ adding a vertex to the embedding and connecting it to each of the $n-1$ vertices can in some cases produce a drawing with the lowest possible number of crossings, and in other cases a different method is required. For example, one can obtain a drawing of $K_5$ with one crossing on the sphere after embedding $K_4,$ whereas a four-crossing drawing of $K_8$ on a torus cannot be obtained by embedding $K_7.$ A four-crossing drawing of $K_9$ on $S_2$ was given by Riskin \cite{riskin}, and this drawing can be obtained from an embedding of $K_8$ by placing the ninth vertex in a quadrilateral face. In this paper, I show that $K_9$ can be embedded on $S_3$ with a hexagon with all vertices distinct, and the other three vertices in faces adjacent to the hexagon. Thus, there is a drawing of $K_{10}$ on $S_3$ with three crossings, and since it was shown by Kainen \cite{kainen} that no smaller number of crossings is possible, I have $cr_3(K_{10})=3$.

\section*{The Heawood Conjecture}
Given a graph $G,$ the \emph{chromatic number} of $G$, denoted $\chi(G)$, is the minimum number of colors needed to color the vertices so that no two neighboring vertices have the same color. In 1890, Percy John Heawood posed a question that would soon spark broad interest in embedding graphs in surfaces: given a surface $S$, consider all graphs $G$ that can be embedded in $S.$ What is the highest chromatic number $\chi(G)$ of all of them? The four color map problem, which asks the same question of the sphere $S_0,$ had been widely known since 1852. Heawood proved that $\chi(G)\leq\lfloor\frac{7+\sqrt{1+48g})}{2}\rfloor$ for all graphs $G$ that can be embedded on a surface of genus $g$, and claimed that for any surface of genus $g$, there is a graph $G$ achieving this upper bound \cite{heawood}. He also showed that the chromatic number of a torus is seven, proving his claim for $g=1$. In 1891, Lothar Heffter described how to embed $K_7$, $K_8$, $K_9$, $K_{10}$, $K_{11}$, and $K_{12}$ in a surface of minimum genus, proving Heawood's claim for $g\leq6$ \cite{heffter}. Gerhard Ringel proved the case $g=7$ by embedding $K_{13}$ in 1952, and he started to consider $K_n$ for the twelve cases of of $n$ modulo $12$. Finally, in 1968, Gerhard Ringel and J. W. T. Youngs proved the Heawood map-coloring theorem by finding embeddings for all residues modulo 12, with the help of Jean Mayer and others for a finite number of cases not covered by their proof \cite{mayer}. Thus Heawood's conjecture was proved for all orientable surfaces. It is also true for non-orientable surfaces except for the Klein bottle. In the course of their work, Ringel and Youngs established a relationship between the genus of a complete graph and chromatic number of a surface. Given $n$, if $p$ is the largest integer such that $\gamma(K_p)\geq n$, then $\frac{(p-3)(p-4)}{12}\geq n\geq\frac{(p-3)(p-2)}{12}$. This inequality is equivalent to $p^{2}-7p+12\geq 12n\geq p^{2}-5p+6$. Solving for $p$, it can be seen that if $\gamma(K_n)=\lceil\frac{(n-3)(n-4)}{12}\rceil$, then $\chi(S_g)=\lfloor\frac{7+\sqrt{1+48g})}{2}\rfloor$ \cite{harary}\cite{ringel}.

\section*{Known Results}
A similar question arose: what is the minimum crossing number of graphs? In 1960, Guy conjectured the minimum crossing number of complete graphs on a genus $0$ surface: $Z(n) = \frac{1}{4}\lfloor\frac{n}{2}\rfloor\lfloor\frac{n-1}{2}\rfloor\lfloor\frac{n-2}{2}\rfloor\lfloor\frac{n-3}{2}\rfloor$ \cite{originalguy} \cite{guy}. The conjecture has been proven for $n\leq12$ \cite{pan}, but the crossing number for $n\geq13$ remains unconfirmed. In 1968, Richard K. Guy, Tom Jenkyns, and Jonathan Schaer proved the lower and upper bounds of the toroidal crossing number of the complete graph as $\frac{23}{210} \binom{n}{4}$ and $\frac{59}{216} \binom{n-1}{4}$ when $n \geq 10$ \cite{guyjenkyns}. Paul C. Kainen proved a lower bound for the crossing numbers of complete graphs, bipartite graphs, and cubical graphs in 1972: $cr_g (K_n)\geq\binom{n}{2}-3n+3(2-2g)$. Kainen conjectured that $cr_g (K_n)$ is equal to this lower bound whenever $g$ is one less than the genus of $K_n$ and $K_n$ does not provide a triangulation of $S_{g+1}$ \cite{kainen}. However, Adrian Riskin disproved Kainen’s conjecture by showing the genus $2$ crossing number of $K_9$ graph is $4$, not $3$ as Kainen’s conjecture suggests \cite{riskin}. In 1995, Farhad Shahrokhi, László A. Székely, Ondrej Sýkora, and Imrich Vrt'o proved a theorem that gives a lower bound for the crossing number. Given a graph $G$, when $\frac{n^2}{e}\geq g$, the crossing number is at least $\frac{ce^3}{n^2}$. When $\frac{n^2}{e}\leq g\leq\frac{3}{64}$, the crossing number is at least $\frac{ce^2}{g+1}$ \cite{10.1007/3-540-57899-4_68}.

\begin{table}[h]
\begin{center}
\begin{tabular}{ |c|c c c c c c c c c c c c| } 
\hline
$n\pmod{12}$ & 0 & 1 & 2 & 3 & 4 & 5 & 6 & 7 & 8 & 9 & 10 & 11 \\
\hline
$f(n)$ & 0 & 3 & 5 & 0 & 0 & 5 & 3 & 0 & 2 & 3 & 3 & 2 \\
\hline
\end{tabular}
\caption{$f(n)$ gives the number of edges that would need to be added to an embedding of $K_n$ on a surface of its minimum genus to obtain a triangulation}
\label{table:1}
\end{center}
\end{table}

\section*{The Genus 3 Case: Program}
Thus any embedding of $K_9$ on $S_3$ has either three quadrilateral faces, a quadrilateral and a pentagon, or a hexagon. Table \ref{table:1} shows the number of edges that would need to be added to obtain a triangulation of the surface for each complete graph $K_n$.

Given an embedding of a graph G with vertices labelled $1,\dots,n$ on the orientable surface $S_g,$ a \emph{rotation sequence} for this embedding consists of $n$ rows, where each row $r_i$ is the sequence of vertices met opposite vertex $i$ along each edge incident to vertex $i$ in clockwise order, starting from a given edge. The Heffter-Edmonds Rotation Principle states that every rotation sequence has a unique orientable graph embedding.

I built an algorithm constructing all rotation sequences for $K_9$ on $S_3$. I determined that there exists an embedding of $K_9$ on $S_3$ that contains a hexagon with $6$ distinct vertices. Furthermore, I know that three vertices are in triangular faces adjacent to this hexagon.

To enumerate all embeddings of $K_9$ on $S_3,$ I generated arbitrary permutations of vertices as rows of candidate rotation sequences and checked:

1. orientability: if an edge shows up in two rows, it is in opposite order;

2. if the number of triangular faces mentioned is less than or equal to $22$;

3. if each edge is in at most two faces;

4. if the faces containing a given vertex make a cycle of length 8.

Heffter already showed an embedding of $K_9$ on $S_3$ with $22$ triangles and one hexagon that contains a vertex mentioned twice. If there is a hexagon with distinct vertices and the other three vertices in faces adjacent to that hexagon, a tenth point inside the hexagon will be able to connect to all six vertices of the hexagon and to the other three vertices with only three crossings. Thus, if I find such an embedding, then the genus 3 crossing number of $K_{10}$ is 3.

Algorithm \ref{algorithm:1} generates all possible rotation sequences of the hexagon and eliminates the impossible drawings. A hexagon with an outer edge showing up in the same order of vertices means that the drawing is not orientable. Sequences with more than $22$ triangular faces or with edges contained in more than two faces can be rejected. Thus, in the end, the program will only print rotation sequences that produce an orientable surface of genus $3$.

Using rotation sequences, I labeled the vertices from $1$ to $9$, giving each vertex a row of different vertices it is connected to in order. In row 1, I may assume that vertices 7, 8, and 9 appear in that order, meaning that there are $\frac{6!}{3!}=120$ permutations to check. In each of the rows 2, 3, 4, 5, and 6, the first vertex after the hexagon must be the same as the last vertex before the hexagon in the previous row. Thus there are $5!=120$ permutations to check for each admissible permutation of the previous row.

This program first generates an array of vertices connected to vertex $1$ excluding those with $3$ as the last element. (Excluding 3 cuts down slightly but is not technically necessary; thus this step is omitted from Algorithm 1.) Then it generates another array of vertices connected to vertex $2$ only if the last element is not $4$ and all outer edges are mentioned in the opposite order compared to those in row $1$. The same procedure is repeated until vertex $4$, when I now check that the number of triangular faces is less than $22$, because Heffter already proved that the number of faces in $K_9$ on $S_3$. For vertex $5$, I check that each edge is mentioned in at most two faces. The same is repeated for vertex $6$. This makes sure that the program does not have to go through generating all sequences but only fully process sequences that have potential to be an embedding. Lastly, I check if the cycle lengths of faces going around vertices 7, 8, and 9 are 8. If so, the program prints the sequences so that I can reproduce them into drawings.

\newtheorem{thm}{Theorem}
\begin{thm}
There are two non-isomorphic embeddings of $K_9$ on $S_3$ with a hexagon with all six vertices distinct.
\begin{proof}
Algorithm \ref{algorithm:1} generated a total of eight rotation sequences with six distinct vertices. Four of them were permutations of the rotation sequence given in Table \ref{table:6}, and four were permutations of the rotation sequence given in Table \ref{table:7}. In Case 2, three faces adjacent to the hexagon had non-hexagon vertices and were adjacent to a triangle with another non-hexagon vertex. This is not true of Case 1.
\end{proof}
\end{thm}

\begin{table}[h!]
\begin{center}
\begin{tabular} { | c | c  c  c  c c c c c | }
 \hline
1&6&3&7&8&5&4&9&2\\
2&1&9&6&4&8&7&5&3\\
3&2&5&9&7&1&6&8&4\\
4&3&8&2&6&7&9&1&5\\
5&4&1&8&9&3&2&7&6\\
6&5&7&4&2&9&8&3&1\\
\hline
\end{tabular}
\caption{First six rows of a rotation sequence embedding $K_9$ on $S_3$, Case 1}
\label{table:6}
\end{center}
\end{table}

\begin{table}[h!]
\begin{center}
\begin{tabular} { | c | c  c  c  c c c c c | }
 \hline
1&6&3&7&8&5&4&9&2\\
2&1&9&8&6&4&7&5&3\\
3&2&5&9&7&1&6&8&4\\
4&3&8&7&2&6&9&1&5\\
5&4&1&8&9&3&2&7&6\\
6&5&7&9&4&2&8&3&1\\
\hline
\end{tabular}
\caption{First six rows of a rotation sequence embedding $K_9$ on $S_3$, Case 2}
\label{table:7}
\end{center}
\end{table}

\begin{algorithm}[h]
\label{algorithm:1}
\newpage 
\caption{Searching through all rotation sequences to find the set of all embeddings of $K_9$ on $S_3$ with a hexagon with $6$ distinct vertices.}
\lstset{numbers=left, numberstyle=\tiny, stepnumber=1, numbersep=5pt}
\SetKwData{Left}{left}\SetKwData{This}{this}\SetKwData{Up}{up}
\SetKwFunction{Union}{Union}\SetKwFunction{FindCompress}{FindCompress}
\SetKwInOut{Input}{input}\SetKwInOut{Output}{output}
\Output{rotation sequences with a hexagon with $6$ distinct vertices}
\BlankLine
check\_perms\_of\_row$(i)$:\\
Let $r_{i,1}=i-1$ and $r_{i,8}=i+1$\\
\ForAll{permutations $p$ of $\{1, 2, 3, 4, 5, 6, 7, 8, 9\}-\{i\}-\{i+1\}-\{i-1\}-\{v\}$}{fill in $r_{i, 3}$ through $r_{i, 7}$ with $p$\\
\tcc{in lines 2 and 3, $i+1$ and $i-1$ are to be read $\mod{6}$ using the numbers $1, 2, 3, 4, 5, 6$}\
	\If{any edge is in two rows in the same direction}{go to next $p$}\
	\If{total number of faces mentioned through row $r_i>22$}{go to next $p$}\
	\If{any edge is in more than two faces}{go to next $p$}\
\If{$i<6$}{check\_perms\_of\_row$(i+1)$}\
\If{$i=6$}{\If{the total number of faces mentioned through row $r_6$ is less than 22}{
add this rotation sequence to set of candidates for further verification}
\If{for any $j \in \{7,8,9\}$, the triangular faces containing $j$ do not make a cycle of length 8}{go to next $p$}
print out the first six rows of the rotation sequence}}
\end{algorithm}

\clearpage
\begin{table}[h]
\begin{center}
\begin{tabular} { | l | r | l}
 \hline
 Check & Number of Cases\\
 \hline
 Row 4: Check if edges opposite & 969598\\
 Check if number of faces $\leq22$ & 261359\\
 \hline
 Row 5: Check if edges opposite & 1169027\\
 Check if number of faces $\leq22$ & 3307\\
 Check if each edge is in at most $2$ faces & 110\\
 \hline
 Row 6: Check if edges opposite & 100\\
 Check if number of faces $\leq22$ & 38\\
 Check if each edge is in at most $2$ faces & 38\\
 \hline
 See if all $8$ edges connecting to $7, 8, 9$ make a cycle & 8\\
\hline
\end{tabular}
\caption{Number of arrangements of the first $i$ rows of candidate rotation sequence for $K_9$ satisfying each criterion}
\label{table:2}
\end{center}
\end{table}

\begin{thm}
The genus $3$ crossing number of $K_{10}$ is $3$.
\begin{proof}
Kainen's inequality $cr_g (K_n)\geq\binom{n}{2}-3n+3(2-2g)$ gives the lower bound of the crossing number of $K_{10}$ on a genus $3$ surface. With Algorithm \ref{algorithm:1}, I found rotation sequences including a hexagon with 6 distinct vertices and the other three vertices in faces adjacent to the hexagon. An additional vertex in the hexagon will be able to connect with the six vertices of the hexagon and cross an edge each for vertices 7, 8, and 9. The drawing is shown in Figure \ref{figure:4}.
\end{proof}
\end{thm}

\section*{The Genus 4 Case: Program}
Although the foundations of the algorithm for $K_{10}$ are similar to those of $K_9$, Algorithm \ref{algorithm:2} is different with one more vertex present. There is, of course, one more place in each row. Most importantly, there are a total of 210 possibilities in the first row, because I may assume vertices 7, 8, 9, and 10 appear in that order and the middle seven elements have to be arranged so that $\frac{7!}{4!}$. Other rows have 720 possibilities, because three elements are known. 

The fact that this algorithm works so plainly for $K_9$ is in part dependent on the fact that no sequences giving fewer than 22 faces among the first six rows turn out to exist. But for $K_{10}$, the analogous algorithm outputs that it is possible to have 26, 27, or 28 faces among the first six rows. I thus divide the algorithm into three cases: when there are 26 faces, 27 faces, and 28 faces. In the case of 28 faces, the algorithm is similar to the algorithm for $K_9$ in that I simply check if the cycle length of each vertex is 9. No configurations existed with this property. For the case of 27 faces, I iterate over triangles consisting only of vertices 7, 8, 9, and 10, adding each to the collection of faces and checking whether it completes an embedding of $K_{10}$. 8 configurations were found to be completable. In the case of 26 faces, I do the same but with pairs of triangles instead of triangles. 11 configurations were found to be completable.

\begin{thm}
The genus $4$ crossing number of $K_{11}$ is $4$.
\begin{proof}
Kainen's inequality $cr_g (K_n)\geq\binom{n}{2}-3n+3(2-2g)$ gives a lower bound of 4 for the crossing number of $K_{11}$ on a genus $4$ surface. With Algorithm \ref{algorithm:2}, I found rotation sequences including a hexagon with 6 distinct vertices and the other four vertices in faces adjacent to the hexagon. An additional vertex in the hexagon will be able to connect with the six vertices of the hexagon and cross an edge each for vertices 7, 8, 9, and 10.

One of the rotation sequences is given in Table \ref{table:8}. The sequence does not mention two faces: faces consisting of vertices 7, 8, and 10 and of vertices 7, 9, and 10. The corresponding embedding is shown in Figure \ref{figure:6}.
\end{proof}
\end{thm}

\begin{table}[h!]
\begin{center}
\begin{tabular} { | c | c  c  c  c c c c c c | }
 \hline
1&6 & 3 &7 &8 &4& 9& 5& 10& 2\\
2&1& 10& 6& 7& 4& 8& 5& 9& 3\\
3&2 &9 &8& 10& 5& 7& 1& 6& 4\\
4&3 &6 &10& 9& 1& 8& 2& 7& 5\\
5&4 &7 &3 &10& 1& 9& 2& 8& 6\\
6&5 &8 &9 &7& 2& 10& 4& 3& 1\\
\hline
\end{tabular}
\caption{First six rows of a rotation sequence embedding $K_{10}$ on $S_4$}
\label{table:8}
\end{center}
\end{table}

\begin{algorithm}[h!]
\label{algorithm:2}
\newpage 
\caption{Searching through all rotation sequences to find the set of all embeddings of $K_{10}$ on $S_4$ with a hexagon with $6$ distinct vertices.}
\lstset{numbers=left, numberstyle=\tiny, stepnumber=1, numbersep=5pt}
\SetKwData{Left}{left}\SetKwData{This}{this}\SetKwData{Up}{up}
\SetKwFunction{Union}{Union}\SetKwFunction{FindCompress}{FindCompress}
\SetKwInOut{Input}{input}\SetKwInOut{Output}{output}
\Output{rotation sequences with a hexagon with $6$ distinct vertices}
\BlankLine
check\_perms\_of\_row$(i)$:\\
Let $r_{i,1}=i-1$ and $r_{i,8}=i+1$\\
\ForAll{permutations $p$ of $\{1, 2, 3, 4, 5, 6, 7, 8, 9, 10\}-\{i\}-\{i+1\}-\{i-1\}-\{v\}$}{fill in $r_{i, 3}$ through $r_{i, 8}$ with $p$\\
\tcc{in lines 2 and 3, $i+1$ and $i-1$ are to be read $\mod{6}$ using the numbers $1, 2, 3, 4, 5, 6$}\
	\If{any edge is in two rows in the same direction}{go to next $p$}\
	\If{total number of faces mentioned through row $r_i>28$}{go to next $p$}\
	\If{any edge is in more than two faces}{go to next $p$}\
\If{$i<6$}{check\_perms\_of\_row$(i+1)$}\
\If{$i=6$}{\If{the total number of faces mentioned through row $r_6$ is less than 28}{
\If{the total number of faces is 27}{check if triangles of vertices 7, 8, 9, 10 can complete an embedding}
\If{the total number of faces is 26}{check if pairs of triangles of vertices 7, 8, 9, 10 can complete an embedding}
\If{the total number of faces is less than 26}{add this rotation sequence to set of candidates for further verification}}
\If{for any $j \in \{7,8,9,10\}$, the triangular faces containing $j$ do not make a cycle of length 9}{go to next $p$}
print out the first six rows of the rotation sequence}}
\end{algorithm}

\begin{table}[h]
\begin{center}
\begin{tabular} { | l | r | l}
 \hline
 Check & Number of Cases\\
 \hline
 Row 4: Check if edges opposite & 868209541\\
 Check if number of faces $\leq28$ & 79975\\
 \hline
 Row 5: Check if edges opposite & 4290114\\
 Check if number of faces $\leq28$ & 3588229\\
 Check if each edge is in at most $2$ faces & 8860\\
 \hline
 Row 6: Check if edges opposite & 187765\\
 Check if number of faces $\leq28$ & 30226\\
 Check if each edge is in at most $2$ faces & 747\\
 \hline
 See if cycle length around $7, 8, 9, 10$ is 9 & 0\\
 27 faces see if completable & 8 \\
 26 faces see if completable & 11 \\
\hline
\end{tabular}
\caption{Number of arrangements of the first $i$ rows of candidate rotation sequence for $K_{10}$ satisfying each criterion}
\label{table:4}
\end{center}
\end{table}

\clearpage
\section*{Summary}
\begin{table}[h]
\begin{center}
\begin{tabular} { | l | c | c | c | c | }
 \hline
 $g$/$n$ & 8 & 9 & 10 & 11 \\
 \hline
 0 & 18 & 36 & 60 & 100  \\
 \hline
 1 & 4 & 9 & 23 & [37, 42] \\
 \hline
 2 & 0 & 4 & [9, 12] & [16, 27] \\
 \hline
 3 & - & 0 & 3 & [10, 14] \\
 \hline
 4 & - & - & 0 & 4 \\
 \hline
 5 & - & - & - & 0 \\
\hline
\end{tabular}
\caption{Range of $cr_g(K_n)$}
\label{table:3}
\end{center}
\end{table}

\begin{figure}[h]
\caption{$K_{10}$ on a surface of genus 2 with 12 crossings}
\centering
\includegraphics[width=12cm]{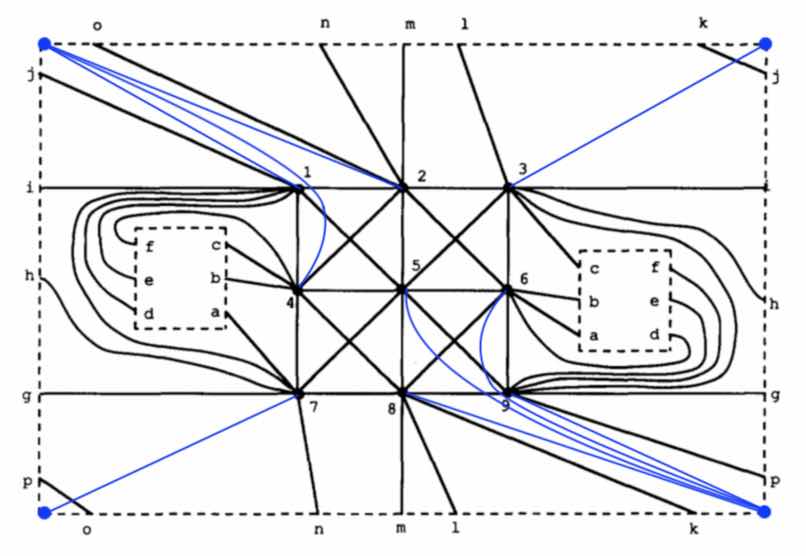}
\label{figure:1}
\end{figure}

\begin{figure}[h]
\caption{$K_{11}$ on a surface of genus 2 with 27 crossings}
\centering
\includegraphics[width=12cm]{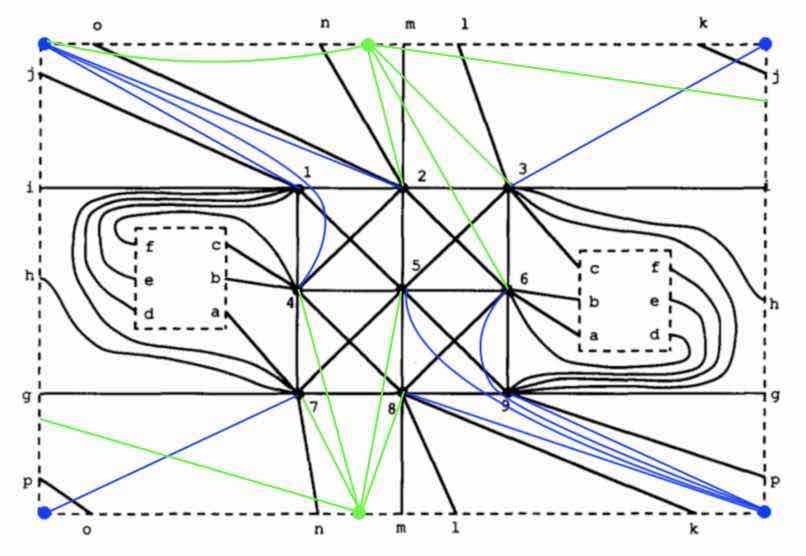}
\label{figure:2}
\end{figure}

The following figures show the drawings of complete graphs on different genus surfaces. Figures \ref{figure:1} and \ref{figure:2} show $K_{10}$ and $K_{11}$ on a genus $2$ surface. The diagrams were based on Adrian Riskin's $cr_2(K_9)=4$ drawing.

\begin{figure}[h]
\caption{Embedding of $K_{9}$ on a surface of genus 3}
\centering
\includegraphics[width=13cm]{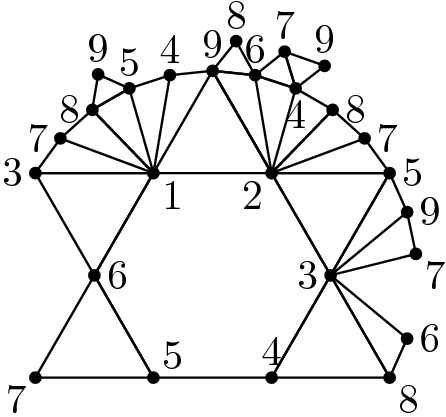}
\label{figure:3}
\end{figure}

\begin{figure}[h]
\caption{$K_{10}$ on a surface of genus 3 with 3 crossings}
\centering
\includegraphics[width=13cm]{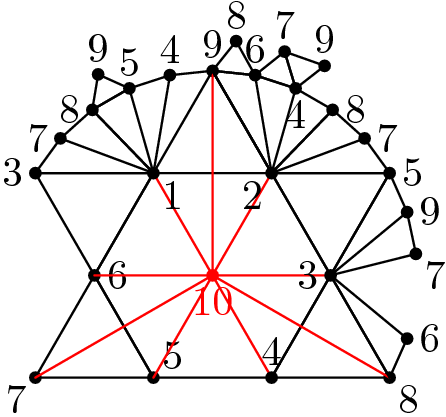}
\label{figure:4}
\end{figure}

\begin{figure}[h]
\caption{$K_{11}$ on a surface of genus 3 with 14 crossings}
\centering
\includegraphics[width=13cm]{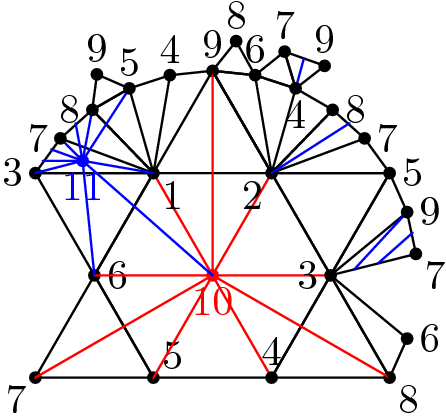}
\label{figure:5}
\end{figure}

Figures \ref{figure:3}, \ref{figure:4}, and \ref{figure:5} show graphs on a genus 3 surface using the rotation sequences. Figure \ref{figure:3} shows an embedding of graph $K_9$ on a surface of genus 3. Figure \ref{figure:4} is graph $K_{10}$ on a genus 3 surface with 3 crossing numbers, which is the lowest crossing number possible. Figure \ref{figure:5} shows $K_{11}$ on a surface of genus 3 with 14 crossings, which I establish as the upper bound of the crossing number. Both Figures \ref{figure:4} and \ref{figure:5} were produced by adding a point on each graph.

\begin{figure}[h]
\caption{Embedding of $K_{10}$ on a surface of genus 4}
\centering
\includegraphics[width=13cm]{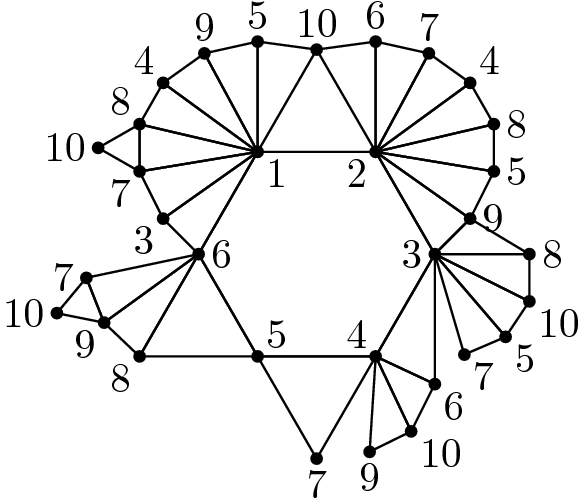}
\label{figure:6}
\end{figure}

\begin{figure}[h]
\caption{$K_{11}$ on a surface of genus 4 with 4 crossings}
\centering
\includegraphics[width=13cm]{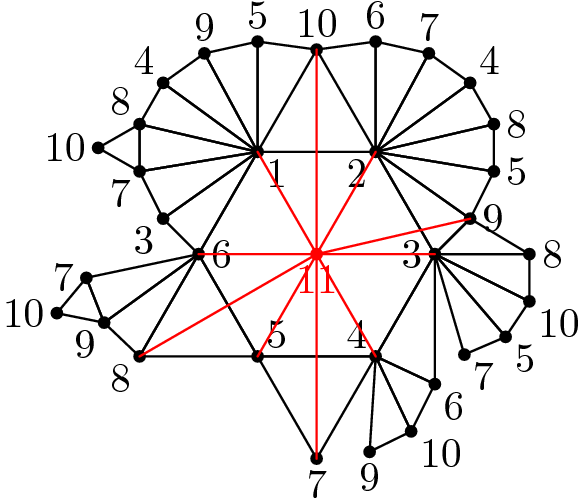}
\label{figure:7}
\end{figure}

Figures \ref{figure:6} and \ref{figure:7} show graphs on a genus 4 surface also using the rotation sequences. Figure \ref{figure:6} is an embedding of $K_{10}$ on a surface of genus 4. Figure \ref{figure:7} shows $cr_4(K_{11})=4$.

\clearpage
The crossing number for each complete graph when $g=0$ has been obtained using Guy's conjecture $Z(n) = \frac{1}{4}\lfloor\frac{n}{2}\rfloor\lfloor\frac{n-1}{2}\rfloor\lfloor\frac{n-2}{2}\rfloor\lfloor\frac{n-3}{2}\rfloor$, which Pan and Richter proved for $n\leq12$. The crossing numbers $4, 9,$ and $23$ for $n=8, 9, 10$ when $g=1$ were proven by Guy, Jenkyns, and Schaer. For $cr_1(K_{11})$, Guy, Jenkyns, and Schaer prove that the upper bound of its crossing number is $42$, while the lower bound is $37$ using $\frac{23}{210} \binom{n}{4}$. The genus of the complete graph was calculated using Ringel and Young's established bound $\gamma(K_n)=\lceil\frac{(n-3)(n-4)}{12}\rceil$, which is why $cr_2(K_8), cr_3(K_9), cr_4(K_{10}),$ and $cr_5(K_{11})$ are all 0. Riskin proved $cr_2(K_9)=4$, and the lower bounds for the rest of the complete graphs on different genus surfaces were established using Kainen's conjecture: $cr_g (K_n)\geq\binom{n}{2}-3n+3(2-2g)$.

The upper bounds for $cr_2(K_{10}), cr_2(K_{11}),$ and $cr_3(K_{11})$ and values of $cr_3(K_{10})$ and $cr_4(K_{11})$ are established above in Theorems 2 and 3 and Figures \ref{figure:1}, \ref{figure:2}, \ref{figure:3}, \ref{figure:4}, and \ref{figure:5}. Table \ref{table:3} shows the lower and upper bounds for $cr_g(K_n)$ when $g\leq5$ and $8\leq n\leq11$.

\clearpage
\bibliographystyle{unsrt}
\bibliography{SampleBib}

\begin{thebibliography}{10}

\bibitem{szekely}
László~A. Székely.
\newblock {Crossing Numbers and Hard Erdős Problems in Discrete Geometry}.
\newblock {\em Combinatorics, Probability and Computing}, 6, 1993.

\bibitem{guy1971}
Richard~K. Guy.
\newblock Latest results on crossing numbers.
\newblock In M.~Capobianco, J.~B. Frechen, and M.~Krolik, editors, {\em Recent
  Trends in Graph Theory}, pages 143--156, Berlin, Heidelberg, 1971. Springer
  Berlin Heidelberg.

\bibitem{heawood}
Percy~John Heawood.
\newblock Map colour theorem.
\newblock {\em Quarterly Journal of Mathematics}, 24, 1890.

\bibitem{riskin}
Adrian Riskin.
\newblock The genus 2 crossing number of $k_9$.
\newblock {\em Discrete Mathematics}, 145, 1995.

\bibitem{kainen}
Paul~C. Kainen.
\newblock A lower bound for crossing numbers of graphs with applications to
  ${K_n}$, ${K_{p,q}}$, and ${Q(d)}$.
\newblock {\em Journal of Combinatorial Theory}, 12, 1972.

\bibitem{heffter}
L.~Heffter.
\newblock Ueber das problem der nachbargebiete.
\newblock {\em Mathematische Annalen}, 38, 1891.

\bibitem{mayer}
Jean Mayer.
\newblock Le problème des régions voisines sur les surfaces closes
  orientables.
\newblock {\em Israel Journal of Mathematics}, 5 1968.

\bibitem{harary}
Frank Harary.
\newblock {\em Graph Theory}.
\newblock Westview Press, 1969.

\bibitem{ringel}
Gerhard Ringel and J.~W.~T. Youngs.
\newblock Solution of the {Heawood Map-Coloring Problem}.
\newblock {\em National Academy of Sciences}, 60, 1968.

\bibitem{originalguy}
Richard~K. Guy.
\newblock A combinatorial problem.
\newblock {\em Nabla (Bulletin of the Malayan Mathematical Society)}, pages
  68--72, 7 1960.

\bibitem{guy}
Richard~K. Guy.
\newblock Crossing numbers of graphs.
\newblock In Y.~Alavi, D.~R. Lick, and A.~T. White, editors, {\em Graph Theory
  and Applications}, pages 111--124, Berlin, Heidelberg, 1972. Springer Berlin
  Heidelberg.

\bibitem{pan}
Shengjun Pan and R~Bruce~Richter.
\newblock The crossing number of ${K_{11}}$ is 100.
\newblock {\em Journal of Graph Theory}, 56:128 -- 134, 10 2007.

\bibitem{guyjenkyns}
Richard~K. Guy, Tom Jenkyns, and Jonathan Schaer.
\newblock The toroidal crossing number of the complete graph.
\newblock {\em Journal of Combinatorial Theory}, 4, 1968.

\bibitem{10.1007/3-540-57899-4_68}
Farhad Shahrokhi, Laszl{\'o}~A. Sz{\'e}kely, Ondrej S{\'y}kora, and Imrich
  Vrt'o.
\newblock Improved bounds for the crossing numbers on surfaces of genus g.
\newblock In Jan van Leeuwen, editor, {\em Graph-Theoretic Concepts in Computer
  Science}, pages 388--395, Berlin, Heidelberg, 1994. Springer Berlin
  Heidelberg.

\end{thebibliography}

\end{document}